\documentclass[a4paper]{article}
\usepackage{amssymb}
\usepackage[dvips]{graphicx}
\usepackage{amsmath}
\usepackage{pstricks}

\newtheorem{remark}{Remark}
\newtheorem{definition}{Definition}

\newtheorem{cor}{Corollary}
\newtheorem{note}{Notation}
\newtheorem{prop}{Proposition}

\title{A time-parallel algorithm for almost integrable
Hamiltonian systems.}
\author{ Hugo Jim\'enez-P\'erez and Jacques Laskar \\ ASD-IMCCE\footnote{ANR Project:---}}
\date{\today}
\begin{document}
\maketitle
\begin{abstract}
  We introduce a time-parallel algorithm for solving numerically
  almost integrable Hamiltonian systems in action-angle
  coordinates. This algorithm is a refinement of that introduced by Saha, 
  Stadel and Tremaine in 1997 (SST97) for the same type of problems. 
  Our refined algorithm has a better convergence obtained from the use
  of derivatives of the perturbing term not considered in the original 
  SST97 algorithm.
  An advantage of this algorithm is its independence of the step-size for 
  the parallelized procedures which can be consider as a particular case
  of the parareal scheme.
\end{abstract}

\section{Introduction}

Many authors agree that the first suggestion of some time-parallel
solution for scalar ordinary differential equation was proposed by
Nievergelt in 1964 and led to the so called multiple-shooting 
methods \cite{Nie1}. 
A few years latter (1967), Miranker and Liniger proposed a family of parallel
Runge Kutta methods for small scale parallelism, based on the 
predictor-corrector method \cite{ML67}. 
 In 1982, 
 Lelarasmee, Ruehli and Sangiovanni-Vicentelli
 introduced 
 the waveform relaxation methods (WR) in 
\cite{LRS82}. 
WR is based 
on the decomposition of a complex system of mixed implicit equations
into a system of single implicit equations. After decomposition, 
each implicit equation can be solved independently of the others and
consequently the system can be solved in parallel in a natural way.
The first implementation of some time-parallel algorithm for ODE
systems which takes advantage of those three methods is 
given by Bellen and Zennaro in \cite{BZ89}. 

In the context of almost integrable Hamiltonian systems, Saha, Stadel 
and Tremaine introduced, in 1997, a parallel method for the 
computation of 
orbits for the Solar System dynamics \cite{SST97}. In part,
their work is the continuation of other papers published in 1992 et 1994
about symplectic integrators and long-term planetary integration 
\cite{ST92, ST94}. 
Simultaneously to Saha \emph{et al.}, Fukushima 
introduced an alternative method for to obtain numerically a global
solution of ODE systems \cite{Fuk1, Fuk2, Fuk3}. 
His method consists in to use the Picard
iteration method to approximate iteratively a global solution. Such 
solutions will be expressed in terms of Chebyshev polynomials to  accelerate 
numerical computations.
Between 1998 and 2000, Erhel and Rault worked on a parallel algorithm
applied to the computation of satellite's trajectories \cite{ER00}.
In the same approach as Bellen and Zennaro, they implement a 
multiple shooting technique but instead of the WR method they use 
the fixed point theory and Newton iterations. They use 
\emph{automatic differentiation} \cite{GRI89} in order to save 
time when computing the Jacobians $J_F\left(u_n^{(k)}\right)$ for 
the Newton iterations.

The parareal algorithm was introduced by Lions, Maday and Turicini in 
\cite{LMT01} and modified by Bal and Maday in \cite{BM02}. It had
received with interest some years latter \cite{AMV09, Bal1, BW08, GV05, 
GH05, Sta1}. This algorithm is based on 
a coarse-discretization predictor (solved sequentially) and a 
fine-discretization corrector (solved in parallel). 
In the very beginning the pararreal algorithm was used to accelerate 
numerical solutions for parabolic and elliptic systems of PDEs.
However, for hyperbolic systems and highly oscillatory problems, the 
coarse solver cannot predict the fine solution in a satisfactory way.
Several attempts has been tested in refinning the coarse solver
\cite{Bal1,CF06,CF08,FC03,AMV09} and the more recent refinement
is a symmetric scheme with projection due to Dai, Le Bris, Legoll and
Maday \cite{DBLM10}. This symmetric algorithm was tested 
for solving almost integrable Hamiltonian systems with good 
results. However, for long-term computations of highly oscillatory 
systems we must to reduce the local error to the size 
$err\sim10^{-16}$, it means, to the $\varepsilon$-machine, and this 
is achieved for several (more than 10) iterations of the parareal
scheme.

In this paper we propose a refinement of the SST97 algorithm 
looking for to accelerate the solution and to preserve the accuracy of
the sequential integrator. In fact, our algorithm converges 100\% to
the sequential underlying integrator althrought the cost of the
corrector step is high. Our current research is about the economy
of the corrector step and the first results are documented in
\cite{JL11}.

 \section{The time-parallel algorithm}
 
The method we introduce in this paper is concerned with numerical solutions
of almost integrable Hamiltonian systems. Although the exposition will be done
for Hamiltonian vector fields with Hamiltonian perturbing part, at this point 
we do not know about any restriction to apply this method to 
vector fields with non-Hamiltonian perturbation.  

This method is an extension of the Saha, Stadel and
Tremaine algorithm \cite{SST97}  and it is based in the multi-shooting and
Picard's iterative methods, as well as the theory of almost integrable 
Hamiltonian systems. We start with the Picard's iterative method for
solving differential systems.

\subsection{Picard's iterative method}

Let's consider the initial value problem (IVP)
\begin{eqnarray}
  \dot y(t)= f(t, y(t)),\qquad y(0)=y_0.
  \label{eqn:sys:1}
\end{eqnarray}
where $y:[0,T]\to\mathbb R^m$ and $f:[0,T]\times \mathbb R^m\to\mathbb R^m$.

Integrating both sides of (\ref{eqn:sys:1})
from zero to $t\in[0,T]$ we obtain the system of differential equations 
written in its integral form
\begin{eqnarray}
  y(t)=y(0) + \int_0^t f(s,y(s))ds, \qquad t\in[0,T].
  \label{eqn:sys:int}
\end{eqnarray}
Applying Picard's iterative method, we can approximate the solution 
for every $t\in[0,T]$ by 
\begin{eqnarray}
  y^{k+1}(t)=y_0 + \int_0^t f(s,y^k(s))ds, \qquad t\in[0,T].
  \label{eqn:sys:iter}
\end{eqnarray}

If $y_0$ belongs to some convex domain $\mathcal D\subset\mathbb R^m$ around 
the solution and if
$F(t,y(t))$ is Lipchitz in $\mathcal D\times [0,T]$, Picard's theorem
assures the convergence of its iterative method. 
An iteration of the Picard's 
method approximates the solution in space and we say that it is a \emph{vertical} 
iteration.

In order to parallelize the numerical computations, we discretize the problem 
by partitioning $[0,T]$ in $N$ small \emph{slices} of size $\Delta t = T/N.$
We set $t_0 = 0$, $t_n = T$ and $t_i=i\Delta t$ such that 

\begin{eqnarray}
   0=t_0<t_1<t_2< \cdots < t_i <\cdots <t_N=T.
\end{eqnarray}
Then the $i$-th
slice is $\Delta t_i=[t_{i}, t_{i+1}]$ for $i=0,\cdots,N-1$.
 In the same way we write $y_n=y(t_n)$ for the value of the solution at time 
 $t_n$ and for the approximations of the Picard's method we will use the 
 superscript $y_n^k= y^k(t_n)$.

Using this discretization we use the linearity of the integral to rewrite  
(\ref{eqn:sys:iter})  like a sum of integrals in the form
\begin{eqnarray}
  y^{k+1}(T) &=& y_0 + \sum_{j=1}^{N} \int_{t_{j-1}}^{t_j} f \left( s,y^k(s) \right)ds, 
  \qquad s\in [t_{j-1},t_j].
  \label{eqn:sys:stadelT}
\end{eqnarray}

The solution at time $t_{n}$ is approximated by the $(k+1)$-th iteration denoted by $y_{n+1}^{k}
=y^{k+1}(t_{n})$ with partial sums 
\begin{eqnarray}
  y^{k+1}(t_n) &=& y_0 + \sum_{j=1}^{n} \int_{t_{j-1}}^{t_j} f \left( s,y^k(s) \right)ds, 
  \quad s\in [t_{j-1},t_j],\quad 1\leq n\leq N.
  \label{eqn:sys:stadeltn}
\end{eqnarray}

Developping the sums we obtain an iterative scheme for the time by
\begin{eqnarray}
  y^{k+1}(t_0) &=& y^k_0 (=y_0)\\
  y^{k+1}(t_{n+1}) &=& y^{k+1}(t_n) + \int_{t_{n}}^{t_{n+1}} f \left( s , y^k(s) \right)ds, 
      \qquad n=0, \dots, N-1.
  \label{eqn:sys:stadel3}
\end{eqnarray}
We call the iterations in time \emph{horizontal} iterations.
The reader must note that, until now, all computations are exacts. The discretization
just give us the way to combine vertical and horizontal iterations. 
\subsection{Perturbed systems}
Since we are interested in almost integrable Hamiltonian systems, we will use perturbed
IVP in the form 
\begin{eqnarray}
   \dot y(t) = f(y(t)) + \epsilon g(t,y(t)), \quad y(0)=y_0,
   \label{eqn:diffsys:feg}
\end{eqnarray}
where $g:[0,T]\times \mathbb R^m\to\mathbb R^m$ is the perturbing function with 
perturbing parameter $\epsilon\ll1$. Its integral expression is
\begin{eqnarray}
    y(t) = y(0) + \int_0^t f(y(s))ds + \epsilon \int_0^t g(s,y(s))ds, \qquad t\in[0,T],
   \label{eqn:diffsys:fegint}
\end{eqnarray}
and the bi-iterative scheme is 
\begin{eqnarray}
  y^{k+1}(t_0) &=& y^k_0 (=y_0)\\
  y^{k+1}(t_{n+1}) &=& y^{k+1}(t_n) + \int_{t_{n}}^{t_{n+1}} f \left( s , y^k(s) \right)
     +	\epsilon g \left( s , y^k(s) \right)ds, 
  \label{eqn:sys:stadel4}
\end{eqnarray}
for $ n=0, \dots, N-1$.
This bi-iterative scheme permit us to compute the definite integral from the preceding iteration 
$k$ and to make a correction with the values of the new iteration $k+1$. It means
that the integrals can be computed in parallel with some numerical scheme and a small
stepsize $\delta t=\Delta t/N_s$, and perform the corrections in sequence, just as in the 
parareal scheme. However, applying this bi-iterative scheme directly will 
introduce a lot of errors to mid and long term computations and it converges very slowly. 
Seen as a \emph{predictor-corrector} scheme, it corresponds to the identity
map as the predictor, which produce a so far approximation.

One way to boost the convergence of the method is to consider the 
problem of solving (\ref{eqn:sys:stadel4}) as a system of algebraic equations. 
We denote the integrals by
\begin{eqnarray*}
  F(y^k) = \int_{t_{n-1}}^{t_n} f(s,y^k(s))ds\\
  G(y^k) = \int_{t_{n-1}}^{t_n} g(s,y^k(s))ds
\end{eqnarray*}
and we consider $a_n:=y^{k+1}(t_n)$ as a fixed parameter (it corresponds to the
initial condition to the sub-problem when $t\in[t_n,t_{n+1}]$). Then, we want to solve 
\begin{eqnarray*}
  y^{k+1} = a_n + F(y^k)+\epsilon G(y^k),
\end{eqnarray*}
as the fixed point problem
\begin{eqnarray}
  y = a_n + F(y) + \epsilon G(y), \qquad y\in \mathbb R^m.
  \label{eqn:fix:point}
\end{eqnarray}

Let $y^*\in \mathbb R^m$ be a solution for (\ref{eqn:fix:point}), 
then for every $y\in\mathbb R^m$ the following identity holds 
\begin{eqnarray*}
   y^* = a_n + F( y^*) + \epsilon G( y) + \epsilon\left( G( y^*)  
   - G(y)\right).
\end{eqnarray*}
It implies that the iteration 
\begin{eqnarray}
  y^{k+1} = a_n + F(y^{k+1}) + \epsilon G( y^k) + \epsilon\left( G( y^{k+1})  
   - G(y^k)\right),
   \label{eqn:ASD1}
\end{eqnarray}
will converges to the solution if $a_n$ is close to $y^*$.
Returning to the integral form for $F$ and $G$ we obtain
the iteration 
\begin{eqnarray}
  y^{k+1}_{n+1} = y^{k+1}_n + \int_{a}^{b} f(s,y^{k+1})ds + \epsilon \int_{a}^{b} g(s, y^k)ds + \epsilon R
   \label{eqn:ASD2}
\end{eqnarray}
where $y^k_n=y^k(t_n)$, $a=t_n$, $b=t_{n+1}$ and
\begin{eqnarray*}
  R = \left( \int_{a}^{b} g(s, y^{k+1})ds 
  - \int_{a}^{b} g(s,y^k)ds\right),
\end{eqnarray*}
is a remainder which goes to zero when $y^k\to y^{k+1}$. 
\begin{remark}
  In \cite{SST97}, the authors considered 
  \begin{eqnarray*}
    \epsilon \left( G(y^{k+1})- G(y^k) \right) \sim 
    \epsilon  DG(y^{k+1})\cdot (y^{k+1}-y^k),
  \end{eqnarray*}
  as a second order term and their bi-iterative scheme is just 
\begin{eqnarray}
  y^{k+1} = a_n + F(y^{k+1}) + \epsilon G( y^k)
   \label{eqn:SST97}
\end{eqnarray}
(expression (9) in \cite{SST97}).  
\end{remark}

The rule (\ref{eqn:ASD1}) define an implicit difference equation which can be 
solved by some iterative scheme for algebraic equations, for example the
WR algorithm \cite{BZ89}. However, for the case of 
almost integrable Hamiltonian systems, there are explicit (symplectic) 
algorithms for which the numerical computation of $F(y^{k+1})$ has no 
additional cost.

\subsection{Almost integrable Hamiltonian systems}
An almost integrable Hamiltonian system in 
generic canonical variables\footnote{Not necessarily \emph{action-angle} 
variables} has a Hamiltonian function in the form 
\begin{eqnarray}
  H(t,p,q) = H_{Int}(p,q) + \epsilon H_{Per}(t,p,q),\qquad p,q\in\mathbb R^r,t\in\mathbb R,
  \label{eqn:Ham:pq}
\end{eqnarray}
where
\begin{itemize}
 \item $H_{Int}(p,q)$ is the integrable term which is independent of the time,
 \item $H_{Per}(t,p,q)$ is the perturbing term which can depends on 
$p$, $q$ and the
time $t$ as well,
  \item $\epsilon$ is the perturbing parameter and it gives the size
of the perturbation, and
  \item $r$ is the number of degrees of freedom (DoF) of the system.
\end{itemize}
 
If we define $y=(p,q)\in\mathbb R^{2r}$ we can write  the Hamiltonian vector field for the almost integrable problem
as
\begin{eqnarray}
  \dot y(t) = J\nabla_yH_{Int}(y(t)) + \epsilon J\nabla_yH_{Per}(t,y(t)), 
         \quad y\in\mathbb R^{2r}, \qquad t \in I\subset\mathbb R,
   \label{eqn:diffsys:Ham}
\end{eqnarray}
where 
\begin{eqnarray}
  J = \left( 
  \begin{array}[h]{cc}
    0 & I \\
    -I & 0
  \end{array}
  \right), \qquad I,0 \in M_{r}(\mathbb R)
\end{eqnarray}
is the canonical symplectic matrix in $M_{2r}(\mathbb R)$ and
$\nabla_y$ means the gradient with respect to $y$. 

From (\ref{eqn:diffsys:Ham}) we see that 
\begin{eqnarray} 
  f(y(t))=J\nabla_yH_{Int}(y(t))\qquad{\rm and} \qquad 
  	g(t,y(t))=J\nabla_y H_{Per}(t,y(t)).
	\label{eqn:f:g}
\end{eqnarray}

Now we return to expression 
(\ref{eqn:ASD2}) and we split it in $p$ and $q$ to obtain two iterative
rules
\begin{eqnarray}
  p^{k+1}_{n+1} &=&  p^{k+1}_n + \int_{a}^{b} f_1(y^{k+1})ds
  	+ \epsilon \int_{a}^{b} 
           g_1(s, y^k)ds + \epsilon R_1
   \label{eqn:ASD3:1}\\
  q^{k+1}_{n+1} &=&  q^{k+1}_n + \int_{a}^{b} f_2(y^{k+1})ds 
  	+ \epsilon \int_{a}^{b} g_2(s, y^k)ds + \epsilon R_2
   \label{eqn:ASD3:2}
\end{eqnarray}
where $a=t_n$, $b=t_{n+1}$, $f=(f_1,f_2)$, $g=(g_1,g_2)$, $R=(R_1,R_2)$ and
\begin{eqnarray*}
  R_i = \left( \int_{a}^{b} g_i(s, y^{k+1})ds 
  - \int_{a}^{b} g_i(s,y^k)ds\right), \qquad i=1,2.
\end{eqnarray*}

The most important fact obtained from these expressions is that 
perturbing terms depend only on the values obtained in the last iteration.
This implies that we can compute them in parallel and attach those values
in the sequential correction.

In the rest of this work we will be interested
in almost integrable systems in action-angle coordinates 
with separable perturbing function in the form
\begin{eqnarray*}
  H(t,p,q) &=& H_{Int}(p) + \epsilon (H^p_{Per}(t,q)+H^q_{Per}(t,p)),
\end{eqnarray*}
in order to use symplectic integrators. 
With this restriction, (\ref{eqn:f:g}) becomes
\begin{eqnarray} 
  f(p)=\frac{\partial H_{Int}}{\partial p}(p),\quad{\rm and} \quad 
  g(t,q,p)=\left(-\frac{\partial H^p_{Per}}{\partial q}(t,q),
      \frac{\partial H^q_{Per}}{\partial p}(t,p)\right).
	\label{eqn:f:g:1}
\end{eqnarray}
Additionally,
we must separate the remainders to be computed in an independent step.
The iterative scheme reduces to 
\begin{eqnarray*}
  p^{k+\frac12}_{n+1} &=&  p^{k+1}_n - \epsilon \int_{a}^{b} 
  	\frac{\partial H^p_{Per}}{\partial q}(s, q^k)ds\\
   q^{k+\frac12}_{n+1} &=&  q^{k+1}_n + \int_{a}^{b} 
   \frac{\partial H_{Int}}{\partial p}(p^{k+\frac12})ds + \epsilon \int_{a}^{b} 
  	\frac{\partial H^q_{Per}}{\partial p}(s, p^k)ds\\
     y_{n+1}^{k+1} &=& y_{n+1}^{k+\frac12} + \epsilon\left( 
     \int_a^b\frac{\partial H_{Per}}{\partial y}(s,y^{k+\frac12})ds -
     \int_a^b\frac{\partial H_{Per}}{\partial y}(s,y^{k-\frac12})ds 
     \right),
\end{eqnarray*}
with the same notation as above.

The two first steps correspond to the SST97 algorithm if we use the 
symplectic implicit integrator of middle point.
We extend the SST97 algorithm using explicit symplectic integrators of
higher order and the second order term $\epsilon R$

We take expression (\ref{eqn:f:g:1}) and we separe $g(t,q,p)$ as
\begin{eqnarray*}
    g_p(t,q)= \frac{\partial H^p_{Per}}{\partial q}(t,q)\qquad {\rm and}
    \qquad g_q(t,p)= \frac{\partial H^q_{Per}}{\partial p}(t,p),
\end{eqnarray*}
where the subscript $_p$ (respectively $_q$) means that $g_p$ (resp. $g_q$)
is the perturbing term of the $\dot p$ (resp. $\dot q$) equation.

Our ``toy'' parallel \emph{leapfrog} scheme is 
\begin{eqnarray*}
    p^{k+\frac12}_{n+\frac{1}{2}} &=&  p^{k+1}_n - \frac{\delta t}{2} \epsilon 
  	g_p(q_n^k)\\
   q^{k+\frac12}_{n+1} &=&  q^{k+1}_n + \delta t 
   \left(f_q(p_{n+\frac12}^{k+\frac12}) +  \epsilon g_q(p_{n+\frac12}^k)\right)\\
   p^{k+\frac12}_{n+1} &=&  p^{k+\frac12}_{n+1} - \frac{\delta t}{2} 
   \epsilon g_p(q_{n+1}^k)\\
    & &------------------\\
     y_{n+1}^{k+1} &=& y_{n+1}^{k+\frac12} + \Delta t\epsilon\left( 
     g(y_{n+1}^{k+\frac12}) -  g(y_{n+1}^{k-\frac12}) \right),
\end{eqnarray*}
which is an explicit scheme. In this case, $\delta t = \Delta t$ 
The reader must note that this scheme is more expensive that the 
original sequential scheme. In order to have a gain on the 
algorithm we must concatenate several \emph{leapfrog} schemes or 
other  explicit symmetric symplectic integrators. 

In fact, we can see this procedure like a parallel composition method 
for symmetric symplectic integrators, equivalent to those used by 
Suzuki \cite{Suz1}, Yoshida \cite{Yos1}, McLachlan \cite{McL1,McL2}
and Laskar-Robutel \cite{LR01}.

For simplicity, we can put together $j$ steps of the 
\emph{leapfrog} scheme to get
\begin{eqnarray*}
    p^{k+\frac12}_{n+\frac{1}{j+1}} &=&  p^{k+1}_n - \frac{\delta t}{2}
        \epsilon g_p(q_n^k)\\
    q^{k+\frac12}_{n+\frac1{j}} &=&  q^{k+1}_n + \delta t 
        \left(f_q\left(p_{n+\frac1{j+1}}^{k+\frac12}\right) +  \epsilon 
    g_q\left(p_{n+\frac1{j+1}}^k\right)\right)\\
    p^{k+\frac12}_{n+\frac{2}{j+1}} &=&  p^{k+\frac12}_{n+\frac{1}
       {j+1}}-\delta t \epsilon g_p\left(q_{n+\frac{1}{j}}^k\right)\\
            \vdots & & \qquad \vdots\\
    q^{k+\frac12}_{n+1} &=&  q^{k+\frac12}_{n+\frac{j-1}{j}} + 
        \delta t \left(f_q\left(p_{n+\frac{j}{j+1}}^{k+\frac12}\right) +  
        \epsilon g_q\left(p_{n+\frac{j}{j+1}}^k\right)\right)\\
    p^{k+\frac12}_{n+1} &=&  p^{k+\frac12}_{n+\frac{j}{j+1}} 
         - \frac{\delta t}{2} \epsilon g_p\left(q_{n+1}^k\right)\\
    & &------------------\\
    y_{n+1}^{k+1} &=& y_{n+1}^{k+\frac12} + \Delta t\epsilon\left( 
    g\left(y_{n+1}^{k+\frac12}\right) -  g\left(y_{n+1}^{k-\frac12}
      \right) \right).
\end{eqnarray*}
In this case we have $\Delta t = j\cdot \delta t$ meanning that
the second order term will be computed once every $j$ $\delta t$-steps.
More complex 
symplectic integrators can be used instead of the \emph{leapfrog}
scheme like the $\mathcal{SABA}_n$ or $\mathcal{SBAB}_n$ integrators 
from Laskar and Robutel \cite{LR01}.

\section{The algorithm}

The original algorithm was designed by Saha \emph{et al.} and introduced 
in \cite{SST97} to compute planetary orbits for Solar System dynamics. The 
integrable part of the Hamiltonian function corresponds to the 
Keplerian orbits of the planets around the Sun (in elliptic coordinates)
and the perturbing term corresponds to the interplanetary forces (in 
Jacobian ccordinates).

Our algorithm uses a high-order multistep integrator where each local
perturbing value must be saved and passed to the master process. In order
to avoid ``bottle-necks'' in the communication it is recommended to
use this algorithm in a shared memory environment. 

\begin{figure}[h]
  \centering
  \includegraphics[scale=0.5]{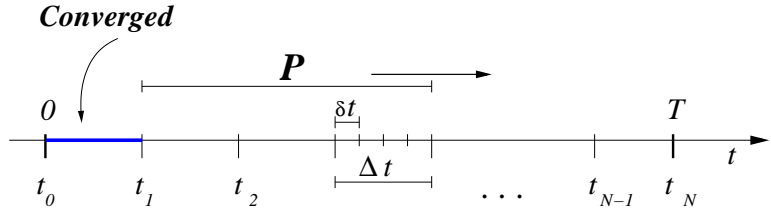}
  \caption{A \emph{ shifting parallelization window} computes $P$ intervals of 
  size 
        $\Delta t$ in parallel. Every interval $\Delta t$ is decomposed
	in smaller intervals of size $\delta t$}
	\label{fig:window}
\end{figure} 

We select an explicit symmetric symplectic integrator $F$, for instance,
the leapfrog or any of the $\mathcal{SABA}_n$ or 
$\mathcal{SBAB}_n$ integrators \cite{LR01}. We concatenate $j>1$ schemes
to be compute in parallel by each thread.
We define $\delta t$ to be the time step of a single $F$ scheme and
$\Delta t=j\cdot \delta t$ the interval computed in parallel. 

Finally, we open a \emph{shifting parallelization window} which contains $P$
intervals of size $\Delta t$. $P$ is the number of threads or, 
if each thread runs on a single processor,
$P$ becomes the number of processors in the parallel machine. 
After each iteration, we verify how many $\Delta t$ intervals 
have converged and we shift the window until the
first non converged $\Delta t$ interval (see Figure \ref{fig:window}).

We denote by 
\begin{eqnarray*}
  \mathfrak F\left(y_{n+\frac ij}^{(k)}\right) &=& c_i\cdot \delta t
  \cdot f_q\left( y_{n+\frac ij}^{(k)} \right) 
\end{eqnarray*}
the $i$-th numerical value for the integrable part and by 
\begin{eqnarray*}
  \mathfrak G\left(y_{n+\frac ij}^{(k)}\right) &=& d_i\cdot \delta t
  \cdot g\left( y_{n+\frac ij}^{(k)} \right) 
\end{eqnarray*}
the corresponding $i$-th value for the perturbing part of the 
Hamiltonian function. The coefficients $c_i$ and $d_i$, $i=1,\dots,j$
are selected based on the symplectic underlying integrator.

We show the general structure of the algorithm in the following 
frame
\begin{center}
  \begin{tabular}[ht]{|l|}
    \hline
    {\bf Algorithm 1.} Jim\'enez-Laskar.\\
    \hline
     {\bf Setup} of the initial guess sequence\\
     $\quad$ $y_0^0=y(0)$, $y_{n+1}^0=\mathfrak F\left(y_n^0\right)$\\
     $\quad$ $G_n^0=0$\\
     {\bf While} $r < N$ {\bf do}\\
     $\quad$ \emph{Parallel} resolution on $[T^n,T^{n+1}]$ 
          for $r\leq n<r+P$: \\ 
     $\qquad$ {\bf For} i=1 {\bf to} j  \\
     $\qquad\qquad$ compute $ 
       \mathfrak F\left(y_{n+\frac ij}^{(k)}\right)
        +\epsilon \mathfrak G \left(y_{n+\frac ij}^{(k)}\right)$,\\ 
     $\qquad\qquad$ save $\mathfrak G\left(y_{n+\frac ij}^{(k)}\right)$,\\
     $\qquad$ {\bf End} (For i),\\
     $\qquad$ $u_{r+1}=y_{r+1}^{(k)}$\\
     $\qquad$ $conv = 1.$\\
     $\qquad$ Converged=TRUE.\\
     $\qquad$ Head=TRUE.\\
     $\quad$ \emph{Sequential} corrections: \\
     $\qquad$ {\bf For} $n=(r+1)$ {\bf to} $(r+P)$\\
     $\quad\qquad$ {\bf For} i=1 {\bf to} j\\
     $\qquad\qquad$ compute $y_{n+1}^{(k+\frac 12)}=  
     \mathfrak F\left( y^{(k+\frac 12)}_{n+\frac ij} \right)
     + \epsilon \mathfrak G \left( y^{(k)}_{n+\frac ij} \right)$.\\
     $\quad\qquad$ {\bf End} (For i)\\
     $\quad\qquad$ compute $G_{n+1}^{(k+\frac 12)} =
     \mathfrak G\left(y_{n+1}^{(k+\frac12)}\right)$\\
     $\quad\qquad$ compute $y_{n+1}^{(k+1)}= y_{n+1}^{(k+\frac12)} + 
       \epsilon \left( G_{n+1}^{(k+\frac12)}
         - G_{n+1}^{(k-\frac12)}\right)$.\\
     $\quad\qquad$ {\bf If} Head=TRUE\\
     $\qquad\qquad$ {\bf If} Converged=TRUE\\
     $\qquad\qquad\quad$ $u_{n+1}= y_{n+1}^{(k+1)}$\\
     $\qquad\qquad\quad$ $conv=conv+1$.\\
     $\qquad\qquad$ {\bf Else}\\
     $\qquad \qquad\quad$ $r = r+conv$ (shifts the window).\\
     $\qquad\qquad\quad$ Converged=FALSE\\
     $\qquad\qquad\quad$ Head=FALSE\\
     $\qquad\qquad\quad$ begins the next parallel iteration\\
     $\qquad\qquad\qquad$ with the available processors.\\
     $\qquad\qquad$ {\bf End} (If Converged=TRUE)\\
     $\qquad\quad$ {\bf End} (If Head=TRUE)\\
     $\qquad$ {\bf End} (For n).\\
     {\bf End} (While).\\
    \hline
  \end{tabular}
\end{center}
where the sequence $\{y_n^{(k)}\}$ is the $k$-th approximation to the solution
$\{u_n\}$. 
It is important to note that this algorithm converges exactly to the 
sequential underlying integrator.

\section{Numerical examples}
We have tested our algorithm with the simple
Hamiltonian pendulum and the Spin-orbit problem. 
We have implemented the code in TRIP \cite{TRIP}
simulating the parallel implementation for several values of $P$
in order to find $P^*$ which optimizes the speed-up function 
$\mathcal S_P$. The discussion about this subject will be treated 
in the next section. 

Hamiltonian problems for very long times. In order to reach this, we have 
tested the algorithm with no tolerance in the error between the sequential 
solution and the parallel aproximation. It means that we test \verb+err==0+
at every iteration where \verb+err+=$\|y_{seq}(t_n)-y^{(k)}(t_n)\|$ 
is the numerical Euclidian distance between the sequential solution 
$y_{seq}(t_n)$ and the $k$-th approximation $y^{(k)}(t_n)$ at time $t_n$. 
Of course, if we relax this condition to the machine tolerance or to the
size of the Hamiltonian error at each $\Delta t$ interval we will have a
gain in the speed-up function.

\begin{figure}[h]
  \begin{center}
    \includegraphics[scale=0.7]{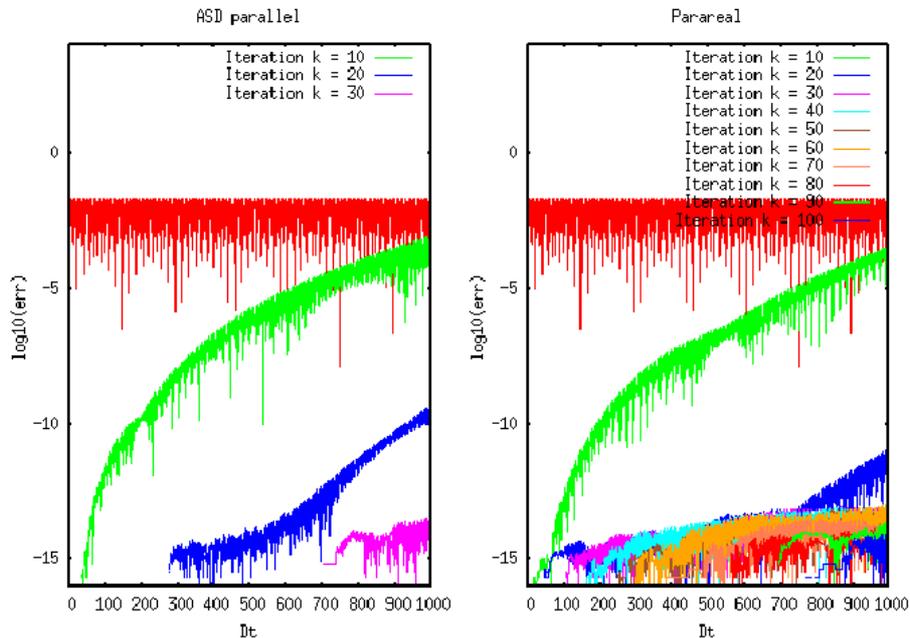}
  \end{center}
  \caption{The Hamiltonian pendulum. 
  Square of the Euclidean distance
  between the sequential and the parallel solutions in 
  $\log_{10}$ scale. Each curve corresponds to the $k$-th iteration for
  $k=10,20,\dots$. The ASD parallel algorithm takes 36 iterations meanwhile 
  the Parareal algorithm takes 109 iterations to total convergence 
  ($err=0$).}
  \label{fig:pend:err}
\end{figure}

\subsection{The simple Hamiltonian pendulum}
The first example is the simple Hamiltonian pendulum with equation
\begin{eqnarray}
  H(p,q) &=& \frac{1}{2}p^2 - \epsilon \cos{q}
  \label{eqn:pend}
\end{eqnarray}
which we integrate
with a $\mathcal{SBAB}_4$ symplectic integrator of order $\mathcal O
(\tau^8\epsilon + \tau^2\epsilon^2)$ \cite{LR01}. 
We have made several test
for different values of $\epsilon$, $\delta t$, $\Delta t$ (in $\delta t$ 
steps), the size of the 
simulation $N$ (in $\delta t$ steps) and the size of the shifting 
parallelization window $P$ ($\Delta t$ intervals in parallel).

\begin{note}
We have de following notation: 
\begin{itemize} 
  \item $C_{\Delta t}$ will denote the mean number of $\Delta t$ intervals 
    converged to the sequential solution at each iteration. 
  \item $I_{(\tau,P)}$ will denote the mean number of iterations needed to make 
    converge $P$ intervals of size $\Delta t$. 
  \item $k_{(\Delta t,P)}$ will denote the total number of iterations 
    needed to obtain the sequential solution.
\end{itemize}
\end{note}

These values are related by the following relationships 
\begin{eqnarray}
  P=C_{\Delta t} \times I_{(\Delta t,P)}\quad{\rm and}\quad 
  k_{(\Delta t, P)}= \frac{N}{jP}I_{(\Delta t,P)} = \frac{N}{jC_{\Delta t}}
  \label{eqn:CIK:ini}
\end{eqnarray}
where $j\in\mathbb N$ is such that $\Delta t = j\delta t$

Table 1 shows values obtained for a simple test with parameters
$\epsilon = 0.01$, $\delta t=0.01$, $\Delta t =1$, $N=10^6$
and the initial conditions $(p_0,q_0)=(1,0)$. In this case, we concatenate
100 $\delta t$ steps in a $\Delta t$ interval and the number $P$ of $\Delta t$ 
intervals computed in parallel varies from 50 to 500 in steps of 50 intervals.

\begin{center}
\begin{tabular}[h]{||l||c|c|c||}
\hline
\hline
P & $C_{\Delta t}$ & $I_{(\Delta t,P)}$ & $k_{(\Delta t,P)}$\\
   \hline
   \hline
        50 & 6.9728 & 7.17072 & 1434\\
        100 & 12.018 & 8.32083 & 832\\
        150 & 16.3918 & 9.15092 & 610\\
        200 & 20.5318 & 9.74097 & 487\\
        250 & 24.3285 & 10.276 & 411\\
        300 & 27.6981 & 10.8311 & 361\\
        350 & 30.6718 & 11.4111 & 326\\
        400 & 33.7804 & 11.8412 & 296\\
        450 & 36.36 & 12.3762 & 275\\
        500 & 38.9066 & 12.8513 & 257\\
   \hline
   \hline
\end{tabular}\\[1em]
Table 1. Values obtained for different sizes of the shifting \\
parallelization window $P$. The values for the parameters are \\
$\epsilon=0.01$, $\delta t=0.01$, $\Delta t=1$, $N=10^6$ and $T=10^4$.

\end{center}

We have made several tests for different values in the parameters and 
we noted that $C_{\Delta t}$, $I_{(\Delta t,P)}$ and $k_{(\Delta t,P)}$ 
have, in all cases, the same  qualitative behavior. 
In particular, $I_{(\Delta t,P)}$ has an almost linear 
behavior then we have fitted a straight line $I(P)\sim I_{(\Delta t,P)}$ 
by least squares technique. Expressions for $C_{\Delta t}$, 
$I_{(\Delta t, P)}$ and $k_{(\Delta t,P)}$ using (\ref{eqn:CIK:ini}) are:
\begin{eqnarray}
  C_{\Delta t}= \frac{1}{j\left( a +\frac{b}{P} \right)},\quad 
  I_{(\Delta t,P)}= j(aP+b)\quad {\rm and}\quad 
   k_{(\Delta t,P)}=\left( a+\frac{b}{P} \right)N.
   \label{eqn:CIK}
\end{eqnarray}

\begin{figure}[h]
  \begin{center}
    \includegraphics[scale=0.55]{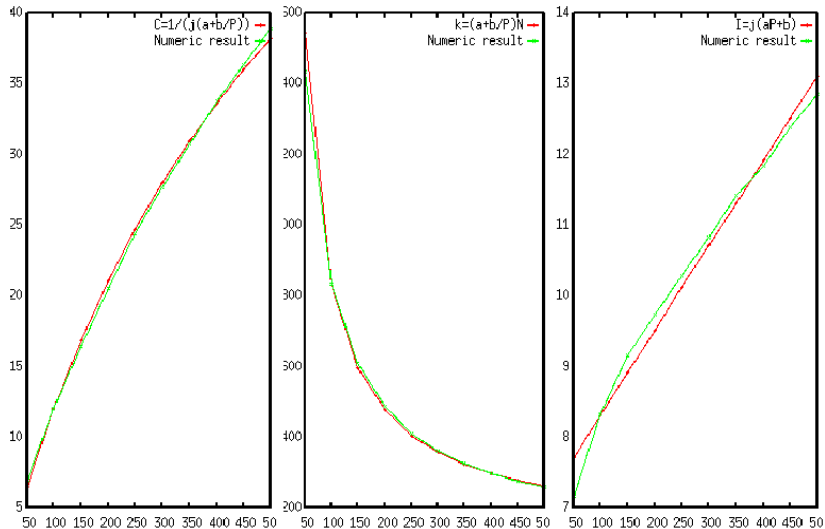}
  \end{center}
  \caption{Plots of the data from Table 1 for $C_{\Delta t}$, $I_{(\Delta t,P)}$,
  and $k_{(\Delta t,P)}$ and its fitting curves with parameters $a=(2.51)^{-4}$
  and $b=(7.64)^{-2}$.}
  \label{fig:Pend}
\end{figure}

Figure \ref{fig:Pend} shows the data of Table 1 with the fitted curves for
the Hamiltonian pendulum with parameters $a=(2.51)^{-4}$ and $b=(7.64)^{-2}$.

\subsection{The spin-orbit problem}
The next example is the spin-orbit problem with equation
\begin{eqnarray}
  H(p,q) &=& \frac{1}{2}p^2-\epsilon\cos{2q}-\alpha(\cos(2q+\phi)-7\cos(2q-\phi))
  \label{eqn:pend2}
\end{eqnarray}
which we integrate with the same symplectic integrator as in the pendulum case.
We have used the parameters $\alpha=10^{-4}$ and $\phi = 0.2$ in 
addition to the values of the last example. 
Table 2 shows values obtained for a simple test of the spin-orbit problem.

\begin{figure}[h]
  \begin{center}
    \includegraphics[scale=0.7]{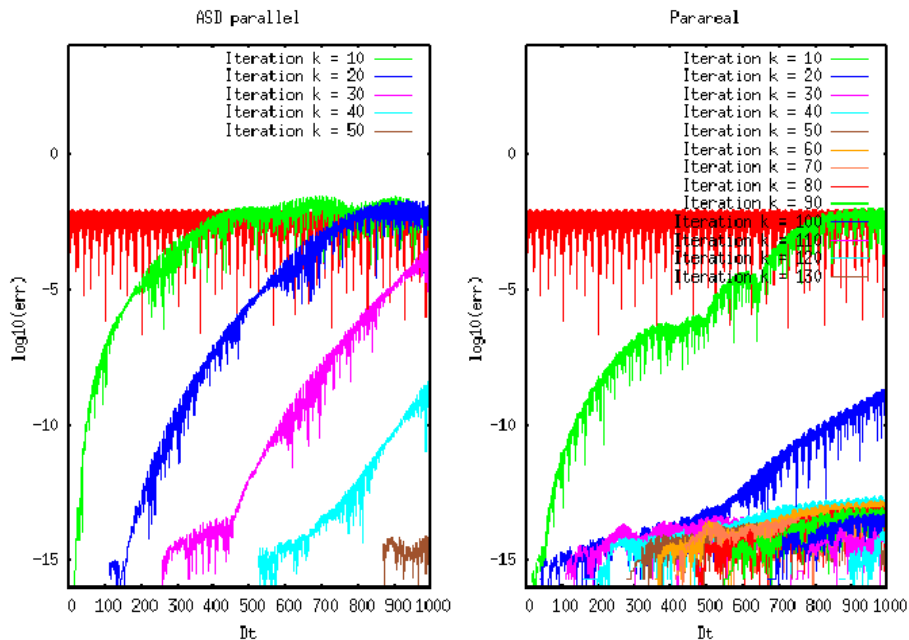}
  \end{center}
  \caption{The spin-orbit problem. Square of the Euclidean distance
  between the sequential and the parallel solutions in 
  $\log_{10}$ scale. Each curve corresponds to the $k$-th iteration for
  $k=10,20,\dots$. The ASD parallel algorithm takes 54 iterations meanwhile 
  the Parareal algorithm takes 135 iterations to total convergence 
  ($err=0$).}
  \label{fig:spin:err}
\end{figure}

\begin{center}
\begin{tabular}[h]{||l||c|c|c||}
\hline
\hline
P & $C_{\Delta t}$ & $I_{(\Delta t,P)}$ & $k_{(\Delta t,P)}$\\
   \hline
   \hline
        50 & 6.08952 & 8.21082 & 1642\\
        100 & 9.81256 & 10.191 & 1019\\
        150 & 12.8028 & 11.7162 & 781\\
        200 & 15.3594 & 13.0213 & 651\\
        250 & 17.6039 & 14.2014 & 568\\
        300 & 19.5293 & 15.3615 & 512\\
        350 & 21.2745 & 16.4516 & 470\\
        400 & 22.9335 & 17.4417 & 436\\
        450 & 23.9211 & 18.8119 & 418\\
        500 & 24.8731 & 20.102 & 402\\
   \hline
   \hline
\end{tabular}\\[1em]
Table 2. Spin-orbit problem: values obtained for different sizes of the\\ 
shifting parallelization window $P$. The values for the parameters are \\
$\epsilon=0.01$, $\alpha=10^{-4}$, $\phi=0.2$, $\delta t=0.01$, $\Delta t=1$, $N=10^6$ and $T=10^4$.
\end{center}

\begin{figure}[h]
  \begin{center}
    \includegraphics[scale=0.6]{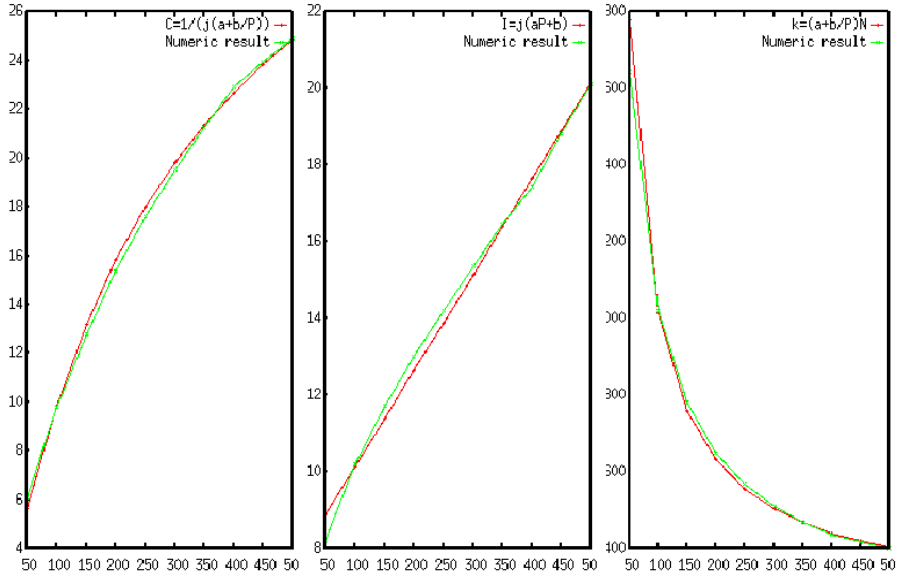}
  \caption{Plots of the data from Table 2 for the spin-orbit problem. 
  We show $C_{\Delta t}$, $I_{(\Delta t,P)}$, and $k_{(\Delta t,P)}$ and
  its fitting curves with parameters
  $a=(1.19)^{-4}$ and $b=(7.11)^{-2}$.}
  \end{center}
  \label{fig:Spin}
\end{figure}

Figure \ref{fig:Spin} shows the data of Table 2 with the fitted curves for
the spin-orbit problem with parameters $a=(1.19)^{-4}$ and $b=(7.11)^{-2}$.

\subsection{Speed-up of the parallel algorithm}
In this subsection we find an expression for the \emph{speed-up}
of the algorithm with parallelization window.
\begin{definition}
  We call \emph{speed-up}, denoted by $S_P(N)$, to the ratio of 
  the time $T_1(N)$ required to solve a given
  problem using the \emph{best known serial} method to the time $T_P(N)$
  required to solve the same problem by a parallel algorithm using $P$ 
  processors
  \begin{eqnarray*}
    \mathcal S_N(P)=\frac{T_1(N)}{T_P(N)}
  \end{eqnarray*}
  where $N$ denotes the size of the problem.
\end{definition}

We use an alternative notation inverting the parameter $P$ and the 
variable $N$ since we are interested in finding some expression 
for the speed-up $\mathcal S_N(P)$ in terms
of $P$. The classical notation is $\mathcal S_P(N)$. 

Let $T_{\mathcal P}$ be the predictor's computing time for a $\Delta t$
interval and let $T_{\mathcal C}$ be the corrector's computing time for
the same interval. We suppose we have a first estimation of the solution.
Then, the speed-up function for the parallel $\mathcal{SABA}_n$ with
shifting parallelization window $P$ is given by
\begin{eqnarray*}
  \mathcal S_N(P) = \frac{NT_{\mathcal P}}{k_{(\Delta t,P)}\left(T_{\mathcal P} 
  	+ PT_{\mathcal C}\right)},
\end{eqnarray*}
where $k_{(\Delta t,P)} \in\mathbb N$ is the number of iterations to converge.
This is the worst case, since with our algorithm we do not wait
for ending the $P$ corrections. In fact, a more accurate expression is
obtained substituting 
$P$ by $(P-C_{\Delta t})$. Since, it is desirable to have an expression 
in terms of $P$ only we use equations (\ref{eqn:CIK}) to obtain such
expression.

\begin{figure}[h]
    \begin{center}
	\includegraphics[scale=0.8]{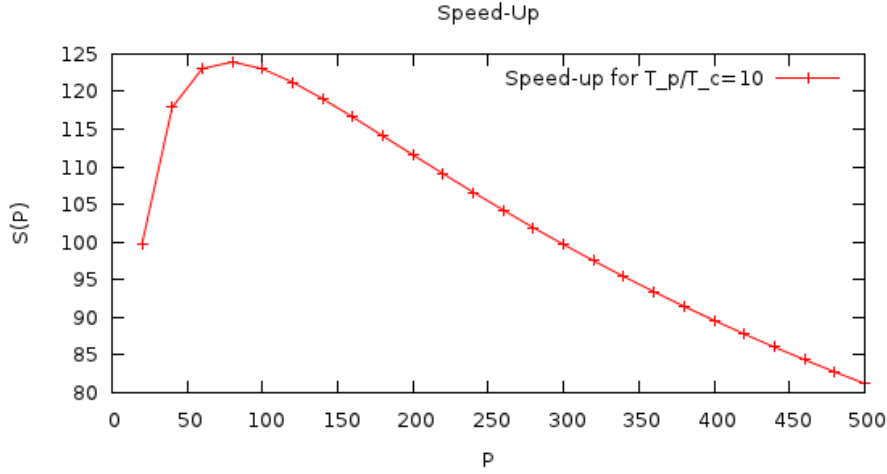}
    \end{center}
    \caption{The Speed-up function with parameters $a=(1.19)^{-4}$, 
    $b=(7.11)^{-2}$, $j=100$ and $\frac{T_{\mathcal P}}{T_{\mathcal C}}=10$.
    }
    \label{fig:SP}
\end{figure}

In what follows, we assume $k_{(\Delta t,P)}=k(P)$ has an analytic expression equivalent to (\ref{eqn:CIK})
as in the cases for the Hamiltonian pendule, the spin-orbit problem 
and the planetary problem (see \cite{JL11}).
After some simplifications and considering $P-C_{\Delta t}$ instead 
of $P$, the formula for the speed-up becomes
\begin{eqnarray}
  \mathcal S_N(P) = \mathcal S(P) = \frac{T_{\mathcal P}}{aT_{\mathcal C}}
   \cdot \frac{P}{ P^2 +B P + C}, \label{eqn:SP}
\end{eqnarray}
where
\begin{eqnarray*}
  B = \left( \frac{T_{\mathcal P}}{T_{\mathcal C}} + \frac{b}{a}
  -\frac{1}{ja}  \right) \quad{\rm and}\quad 
   C= \frac{bT_{\mathcal P}}{aT_{\mathcal C}}.
\end{eqnarray*}
Expression (\ref{eqn:SP}) does not depend on the size of the problem $N$
(see Figure \ref{fig:SP}).

Since we have a quadratic polynomial in the denominator, we may have 
vertical asymptotes in its roots. We impose the condition $\mathcal S(P)>0$
for $P>0$ which implies $P^2+BP+C>0$. 
Additionaly, we have that $\mathcal S(0)=0$, $S(P)>0$ if $P>0$ and since 
$\lim_{P\to\infty}\mathcal S(P)=0$ then there exists a critical point $P^*$
which maximize $\mathcal S(P)$.
We have the following
\begin{prop}
  If the function $k(P)$ has the form (\ref{eqn:CIK}), then $\mathcal S(P)$
  has a maximum in $P^*\in\mathbb N$ which optimize the parallel
  algorithm. Moreover, the theoretic optimal speed-up is 
  \begin{eqnarray}
    \mathcal S(P^*) \cong \frac{1}{ a + 2a\sqrt{\frac{bT_{\mathcal C}}
    {aT_{\mathcal P}}} +\left( b-\frac{1}{j} \right)
    \frac{T_{\mathcal C}}{T_{\mathcal P}}}
    \label{eqn:S:max}
  \end{eqnarray}
\end{prop}
{\it Proof.} We consider expression (\ref{eqn:SP}) 
in the extended case $\mathcal S:\mathbb R \to \mathbb R$ and we 
obtain a differentiable function for $P>0$ such that the original problem 
is just
the restriction $\mathcal S|_{\mathbb N}$. We procede in the classical way 
looking for the critical points by differentiation. Since
\begin{eqnarray*}
   \frac{d\mathcal S}{dP}(P) &=& 
      \frac{ T_{\mathcal P}}{a T_{\mathcal C}} 
       \frac{C-P^2}{(P^2 + BP +C)^2},
\end{eqnarray*}
Then the only positive critical (in fact the maximum) point for the extended
problem is $P_*=+\sqrt{\frac{bT_{\mathcal P}}{aT_{\mathcal C}}}$.
We define 
\begin{eqnarray*}
    P^* &=& \left\{ 
    \begin{array}{ll}
      \left[P_*\right] & {\rm if}\hspace{20pt}
      		\mathcal S_N([P_*]) -\mathcal S([P_*]+1)\leq0\\
		\left[P_*\right]+1 & \hspace{25pt}{\rm otherwise}.
    \end{array}
    \right.
\end{eqnarray*}
where $[\cdot]$ is the maximum integer function.
Finally, for the case where $\mathcal S([P_*])=\mathcal S([P_*]+1)$ we know
that $[P_*]$ is optimal then $P^*$ is well defined.

The second affirmation is obtained directly by substituting $P^*=
\sqrt{\frac{bT_{\mathcal P}}{aT_{\mathcal C}}}$ directly in (\ref{eqn:SP}).

$\hfill\square$

Once we know the optimal value for $P^*$ as a function of $\frac{T_{\mathcal P}}
{T_{\mathcal C}}$, we are interested in the gain, the speed-up which
we can reach with this strategy of parallelizing windows.
We write the speed-up function (\ref{eqn:S:max}) as a function of
the ratio $\frac{T_{\mathcal P}}{T_{\mathcal C}}$ with parameters $a$, $b$
and $j$ as 
\begin{eqnarray}
  \mathcal S_{opt}(P^*) &=& \frac{1}{a + 2\frac{b}{P^*} + \left( b- \frac{1}{j}
  \right) \frac{b}{a(P^*)^2}},
  \label{eqn:S:opt}
\end{eqnarray}

From expression (\ref{eqn:S:opt}) we obtain another useful result
\begin{cor}
  Let $k_{(\Delta t,P)}$ be the number of iterations needed to convergence
  for the parallel algorithm where $k_{(\Delta t,P)}$ is as in
  (\ref{eqn:CIK}) and the parameters $a$, $b$ and $j$ are fixed for some
  particular problem. The speed-up of the parallel algorithm
  is bounded by
  \begin{eqnarray}
    \mathcal S(P) < \frac{1}{a}.
    \label{eqn:S:a}
  \end{eqnarray}
  \label{cor:1}
\end{cor}
{\it Proof.} Take the limit for expression (\ref{eqn:S:opt}) when 
$P^*\to\infty$.

$\hfill\square$

\section{Additional discussion}
What we gain with this algorithm is the time computation of the perturbing
terms. 
Since we have concatenated 100 $\mathcal {SBAB}_4$ schemes, we have called 
500 times the perturbing function in each parallel step (5 times in each 
$\mathcal {SBAB}_4$). In the sequential correction we gain in time 
$500 C_{\Delta t} T_{\mathcal B}$ where $T_{\mathcal B}$ is the time
used in the computation of the perturbing term. However, we must call
the integrable function $400P$ times in the correction step. It means
that this algorithm works in problems where $T_{\mathcal A} \ll
T_{\mathcal B}$ like the spin-orbit problem. 
This is clear from the optimal value of the shifting window 
$P^*=\sqrt{\frac{aT_{\mathcal P}}{bT_{\mathcal C}}}$ in expression 
(\ref{eqn:S:opt}) and the fact that 
$T_{\mathcal P}=j(4T_{\mathcal A}+5T_{\mathcal B})$ and 
$T_{\mathcal C} = 4jT_{\mathcal A}+T_{\mathcal B}$. We have
\begin{eqnarray*}
    \frac{T_{\mathcal P}}{T_{\mathcal C}}= 
    \frac{j(4T_{\mathcal A}+5T_{\mathcal B})}{4jT_{\mathcal A}+T_{\mathcal B}}
    = 1 + \frac{(5j-1)T_{\mathcal B}}{4jT_{\mathcal A}+T_{\mathcal B}} 
    \sim 1 + \frac{T_{\mathcal B}}{T_{\mathcal A}}.
\end{eqnarray*}

Moreover, the time of one 
iteration for this algorithm using the $\mathcal{SBAB}_n$ scheme 
as underlying integrator is 
\begin{eqnarray*}
    T_{\mathcal P} + (P-C_{\Delta t})T_{\mathcal C} &=&  
    j(4T_{\mathcal A} + 5T_{\mathcal B}) + (P-C_{\Delta t})( 4jT_{\mathcal A}  
  + T_{\mathcal B})\\
  &=&  4j(P -C_{\Delta t} + 1)T_{\mathcal A} + (4j+P-C_{\Delta t}) T_{\mathcal B}
\end{eqnarray*}
where $j$ and $P$ must be selected such that the following conditions 
are fulfilled:
\begin{enumerate}
    \item $T_{\mathcal P}\cong (P-C_{\Delta t})T_{\mathcal C}$, this is 
	important to avoid, at maximum, the \emph{dead-time} in the threads. 
    \item $T_{\mathcal P}+(P-C_{\Delta t})T_{\mathcal C}$ be minimal, in 
	order to obtain a large speed-up.
\end{enumerate}

The reader must note that for large values of $P$ or $j$ the 
algorithm might not work if $T_{\mathcal A}\sim T_{\mathcal B}$ 
are comparable. In those cases we have used an alternative corrector 
scheme to parallelize the La2010 solution (see \cite{L11}),
such scheme reduces the time $T_{\mathcal C}$ but increases the
number of iterations needed to convergence (see \cite{JL11}).
Finally, for higher dimensional
Hamiltonian problems as the Solar system dynamics, we find that 
the value of $j$ is very restricted. Some tests with the La2010 
solutions \cite{L11} accept only a maximum of $8$ 
$\mathcal{SABA}_4$ schemes concatenated in a $\Delta t$ interval, 
and for $j>8$ the value of $C_{\Delta t}\equiv1$ which is a worst 
case than the sequential solution.

\end{document}